\documentclass[12pt]{amsart}
\usepackage{amsmath,amsthm,amsfonts,amssymb}
\begin{document} 
\newcommand{\B}{{\mathbb B}}
\newcommand{\C}{{\mathbb C}}
\newcommand{\N}{{\mathbb N}}
\newcommand{\Q}{{\mathbb Q}}
\newcommand{\Z}{{\mathbb Z}}
\renewcommand{\P}{{\mathbb P}}
\newcommand{\R}{{\mathbb R}}
\newcommand{\rc}{\subset}
\newcommand{\rank}{\mathop{rank}}
\newcommand{\trace}{\mathop{tr}}
\newcommand{\dimc}{\mathop{dim}_{\C}}
\newcommand{\Lie}{\mathop{Lie}}
\newcommand{\Auto}{\mathop{{\rm Aut}_{\mathcal O}}}
\newcommand{\alg}[1]{{\mathbf #1}}
\newtheorem*{definition}{Definition}
\newtheorem*{claim}{Claim}
\newtheorem{corollary}{Corollary}
\newtheorem*{Conjecture}{Conjecture}
\newtheorem*{SpecAss}{Special Assumptions}
\newtheorem{example}{Example}
\newtheorem*{remark}{Remark}
\newtheorem*{observation}{Observation}
\newtheorem*{fact}{Fact}
\newtheorem*{remarks}{Remarks}
\newtheorem{lemma}{Lemma}
\newtheorem{proposition}{Proposition}
\newtheorem{theorem}{Theorem}
\title{%
A Lie Group without universal covering
}
\author {J\"org Winkelmann}
\begin{abstract}
We present an example of a disconnected Lie group for which there
is no universal covering (as Lie group).
\end{abstract}
%
\address{%
J\"org Winkelmann \\
 Institut Elie Cartan (Math\'ematiques)\\
 Universit\'e Henri Poincar\'e Nancy 1\\
 B.P. 239\\
 F-54506 Vand\oe uvre-les-Nancy Cedex\\
 France
}
\email{jwinkel@member.ams.org\newline\indent{\itshape Webpage: }%
http://www.math.unibas.ch/\~{ }winkel/
}
\maketitle

A connected and locally pathwise connected topological space $X$ with
a base point $x$ admits a {\em universal covering} $\tilde X$ 
which is constructed
as space of equivalence classes of curves starting at the given base point.
The constant curve with value $x$ defines a base point $\tilde x\in\tilde X$.
Continuous maps between two such spaces mapping the base point of the
source space to the base point of the target space lift
uniquely to a continuous map between the respective universal coverings
mapping basepoint to basepoint.

Using this, for  every connected Lie group the universal covering
admits a structure as Lie group, because the maps defining the
group structure lift uniquely.
The neutral element of the Lie group is the canonical choice for the
basepoint.

However, the pictures changes if we regard disconnected spaces.
In order to construct a universal covering for a disconnected
locally pathwise connected topological space we need a basepoint
for each connected component, because we need to treat each
connected component separately. But for a Lie group only the
connected component containing the neutral element $e$ has a 
canonical choice of a basepoint, namely $e$. For the other
connected components there is no natural choice.
As a consequence, we may construct a universal covering as a manifold
for a disconnected Lie group, but there is no natural way to put
a group structure on this manifold, because there is no canonical
way to lift the map defining the group structure.

In fact, there is an example, where there does not exist any compatible
group structure.
\begin{proposition}
Let $\tilde X=\Z_2\times\Z_2\times\R$ 
and $X=\Z_2\times\Z_2\times(\R/\Z)$
as manifolds and $\pi:\tilde X\to X$ the natural projection
induced by the natural map
$\R\to\R/\Z$. (Here $\Z_2=\Z/2\Z$.)

Let $X$ be endowed with the following group structure:
\[
(a,b,v)\cdot(c,d,w)=(a+c,b+d,v+w+\frac{1}{2}(ad))
\]
with $a,b,c,d\in \Z_2$ and $v,w\in\R/\Z$.

Then there does not exist any Lie group structure on $\tilde X$
for which $\pi$ is a group homomorphism.
\end{proposition}
\begin{proof}
Assume the contrary.
We begin by observing that the connected component $X^0$ 
of $X$ containing $e$ is central in $X$. 
It follows that $ghg^{-1}h^{-1}\in\pi^{-1}(e)$
for all $g\in\tilde X$, $h\in(\tilde X)^0$.
Since $\zeta:(g,h)\to ghg^{-1}h^{-1}$ is continuous, $(\tilde X)^0$
is connected, $\pi^{-1}(e)$ is discrete and $\zeta(g,e)=e$ it follows
that $\zeta(g,h)=e$ for all $g\in\tilde X$, $h\in(\tilde X)^0$.
In other words, $(\tilde X)^0$ is central.

Now let $g\in\tilde X\setminus(\tilde X)^0$. Then $g^2\in (\tilde X)^0$,
because $\tilde X/(\tilde X)^0\simeq X/X^0\simeq \Z_2\times\Z_2$.
Since $(\tilde X)^0$ is central and $(\tilde X)^0\simeq(\R,+)$,
it follows that every connected component of $\tilde X$ contains
a unique element of order two.

Let $g\in\tilde X$. Then conjugation by $g$ stabilizes each of the
connected components, because $\tilde X/(\tilde X)^0$ is commutative.
Therefore conjugation by $g$ must stabilize the unique element of
order two for every connected component. 
Furthermore, conjugation by $g$ acts trivially on the
central subgroup $(\tilde X)^0$.
The elements of order two together $(\tilde X)^0$
generate the whole group $\tilde X$.
Hence conjugation by $g$ is entirely trivial. Since $g$ was arbitrary,
we have deduced that $\tilde X$ is commutative.
But this is in contradiction to the fact that $X$ is not commutative.
\end{proof}

\end{document}